\documentclass[draft]{amsart} 

\usepackage[english]{babel}
\usepackage{xspace}

\newtheorem{theorem}{Theorem}[section]
\newtheorem{proposition}[theorem]{Proposition}

\theoremstyle{definition}
\newtheorem{definition}[theorem]{Definition}
\theoremstyle{remark}

\newtheorem{example}[theorem]{Example}

\newcommand{\Cc}{{\mathbf C}}
\newcommand{\setof}[2]{\left\{\,#1\mid #2\,\right\}}
\newcommand{\Gg}{{\mathbf G}}
\newcommand{\Set}{\ensuremath{\mathbf{Set}}\xspace}
\newcommand{\SimpSets}[1][\Set]{\ensuremath{s #1}\xspace}
\newcommand{\SCat}{\ensuremath{\mathbf{SCat}}\xspace}
\newcommand{\Dd}{{\mathbf D}}
\newcommand{\Ee}{{\mathbf E}}

\newcommand{\Dcat}{\ensuremath{\boldsymbol{\Delta}}\xspace}
\newcommand{\Topo}{\ensuremath{\mathbf{Top}}\xspace}
\newcommand{\op}{\mathrm{op}}
\newcommand{\induc}[3]{{#1}\!\uparrow_{#2}^{#3}}
\newcommand{\commac}[3]{#1\mathord{\backslash}#2\mkern-1.0mu\mathord{\slash}#3}
\newcommand{\comma}[2][F]{#2\mathord{\backslash}#1}
\newcommand{\commaop}[2][F]{#1\mkern-1.0mu\mathord{\slash}#2}

\renewcommand{\(}{\bigl(}
\renewcommand{\)}{\bigr)}

\DeclareMathOperator{\hocolim}{hocolim}
\DeclareMathOperator{\colim}{colim}
\DeclareMathOperator{\obj}{obj}

\DeclareMathOperator{\coend}{coend}
\DeclareMathOperator{\diag}{diag}

\title[Action of a Group on a Diagram]{The Action
  by Natural Transformations of a Group on a Diagram of
  Spaces}

\date{\today}

\author{Rafael Villarroel-Flores}

\address{Instituto de Matem\'aticas, UNAM (Unidad Cuernavaca)\\
Av.~Universidad s/n, col. Lomas de Chamilpa \\
C.~P.~62210 Cuernavaca, Morelos\\
MEXICO}


\begin{document}

\begin{abstract}
  For $\Cc$ a $G$-category, we give a condition on a diagram of
  simplicial sets indexed on $\Cc$ that allows us to define a natural
  $G$-action on its homotopy colimit, and in some other simplicial
  sets and categories defined in terms of the diagram. Well-known
  theorems on homeomorphisms and homotopy equivalences are generalized
  to an equivariant version. 
\end{abstract}

\maketitle

\section{Introduction}
\label{sec:1}

Let $G$ be a group, and $\Cc$ be any small category. Consider a
$\Cc$-diagram of simplicial sets, where the values of the diagram have
a $G$-action.  Then several structures defined in terms of the
diagram, like the colimit and the homotopy colimit, have an induced
structure of $G$-object.  However, it is often the case that one has a
diagram $F\colon \Cc\to\Dd$ where $\Cc$ is a small $G$-category, $\Dd$
is an arbitrary category and the values of $F$ do not necessarily have
a $G$-action, however the homotopy colimit of $F$ does have it.  This
situation was considered in~\cite{MR2002h:55017}, and independently,
by this author in his Ph.~D. thesis~\cite{myphdthesis}, where the
concept of an action of a group $G$ on a functor $F$ by natural
transformations is introduced.  Here we define it formally in
section~\ref{sec:definition-examples}, after the basic definitions in
section~\ref{sec:intoduction}. 

We show that there are induced $G$-actions on colimits, coends, and
bar and Grothendieck constructions of functors on which $G$ acts by
natural transformations.  In section~\ref{sec:homotopy-colimit} we
consider the homotopy colimit, and show some basic identities
involving the constructions defined so far. In section~\ref{sec:2} we
prove the equivariant homotopy invariance of the bar construction.
Finally, in section~\ref{sec:3} we prove the equivariant versions of
the four theorems listed in~\cite[p.~154]{MR86c:55010a} about the
homotopy colimit. Some of them were noted in~\cite{MR2002h:55017},
however a mild additional hypothesis lets us obtain a more precise
result. 

Some of the proofs are as those in~\cite{MR93e:55015}, adapted for the
case of the group action. However we include more details than in the
cited paper, given that homotopy colimit methods have recently been
used by non-topologists, see for example~\cite{wzz}. 

\section{Preliminaries}
\label{sec:intoduction}

Let $G$ be a finite group. We will denote by $\Gg$ the category with a
single object~$*$, in which $\hom_\Gg(*,*)=G$ and the composition
corresponds to group multiplication. For $n\ge 0$, let
$[n]$ be the category associated to the poset~$\{0,1,\ldots,n\}$ with
the usual order.

Let $\SCat$ be the category of small categories and $\Dcat$ be the
full subcategory of~$\SCat$ with objects $\obj\Dcat=\setof{[n]}{n\ge
  0}$. If $\Cc$ and $\Dd$ are categories, with $\Cc$ small, we denote
by $\Dd^\Cc$ the category of functors $\Cc\to\Dd$ (\cite[page
40]{categories}). The category of simplicial sets (see \cite{may}),
denoted $\SimpSets$, is equal to $\Set^{\Dcat^\op}$. The category of
small $G$-categories is defined as $\SCat^{\Gg}$. We identify a small
$G$-category with the image of the functor $\Gg\to\SCat$. From now on,
$\Cc$ will denote a small $G$-category. Note that for each $g\in G$ we
have a functor $g\colon\Cc\to\Cc$, and composition of functors
correspond to group multiplication.  We consider the nerve functor
$N\colon\SCat\to\SimpSets$, given by
$\Cc\mapsto(\hom_{\SCat}(-,\Cc)\colon {\Dcat^\op}\to\Set)$. The nerve
functor sends $G$-categories to $G$-simplicial sets. There is also a
\emph{geometric realization} functor $|\cdot|\colon
\SimpSets\to\Topo$, that sends $G$-simplicial sets to $G$-topological
spaces. We denote $|N(\Cc)|$ simply as $|\Cc|$. 

If $\Dd$ is any category, an object $D$ in it is called a
\emph{$G$-object} if there is a collection of $\Dd$-maps $\{g\colon
D\to D\}$, indexed by the elements of $G$, such that the map
corresponding to the identity element is the identity, and composition
of maps corresponds to group multiplication. 

If $X$ and $Y$ are $G$-topological spaces, a \emph{$G$-homotopy} from
$X$ to $Y$ is a continuous map $H\colon X\times [0,1]\to Y$ such that
$H(gx,t)=gH(x,t)$ for all $g\in G$, $x\in X$ and $t\in [0,1]$. Two
$G$-maps $f_1,f_2\colon X\to Y$ are $G$-homotopic if there is a
$G$-homotopy $H$ from $X$ to $Y$ such that $H(x,0)=f_1(x)$ and
$H(x,1)=f_2(x)$. In this case we write $f_1\simeq_G f_2$.  The
$G$-topological spaces $X$ and $Y$ are $G$-homotopy equivalent if
there are $G$-maps $f\colon X\to Y$ and $f'\colon Y\to X$ such that
$f'f\simeq_G 1_X$ and $ff'\simeq_G 1_Y$. 

We say that two $G$-categories $\Cc_{1}$, $\Cc_{2}$ are
$G$-homotopy equivalent if the spaces $|\Cc_{1}|$, $|\Cc_{2}|$ are. It
is not required that the map $|\Cc_{1}|\to |\Cc_{2}|$ defining the
homotopy equivalence is induced from a functor $\Cc_{1}\to\Cc_{2}$.

If we have a functor $F\colon\Cc\to\Dd$ together with natural
transformations $\{\eta_{g}\colon F\to Fg\}_{g\in G}$ such that
$\eta_{1}$ is the identity and $\eta_{gg'}=\eta_{g}\eta_{g'}$,
Jackowski and S{\l}omi\'nska call $F$ a \emph{right $G$-functor}
\cite{MR2002h:55017}.  Independently, I defined and used the same
concept in my Ph.~D.~thesis \cite{myphdthesis}, and said that in such
situation, $G$ \emph{acts on $F$ by natural transformations}, or
simply that $G$ \emph{acts on} the functor $F$.
% Both terminologies are somewhat equivalent, since, for example, one
% calls a set on which $G$ acts a (right) $G$-set. 
In this paper, we will use both terms
indistinctly. 

In the case that $\Cc$ and $\Dd$ are small $G$-categories, and
$F\colon \Cc\to\Dd$ is a functor such that $F(gC)=gF(C)$,
$F(g\phi)=gF(\phi)$ for all $g\in G$, $C\in\obj\Cc$ and all
$\Cc$-morphisms $\phi$, we will say that $F$ is an \emph{equivariant
functor}.

\section{Definition and Examples}
\label{sec:definition-examples}

We define now our main subject of study in detail:

\begin{definition}
  \label{definition:1}
  Let $F\colon\Cc\to\Dd$ a functor, where $\Cc$ is a small
  $G$-category and $\Dd$ is an arbitrary category.  Suppose that we
  are given a family of $\Dd$-maps $\eta=\{\eta_{g,X}\colon F(X)\to
  F(gX)\}$ indexed by $g\in G$ and $X\in\obj\Cc$ such that
  \begin{enumerate}
  \item $\eta_{1,X}=1_{F(X)}$ for all $X\in\obj\Cc$,
  \item $\eta_{g_1,g_2X}\circ\eta_{g_2,X}=\eta_{g_1g_2,X}$
    for any $X\in\obj\Cc$, $g_1,g_2\in G$,
  \item $\eta_{g,Y}\circ F(f)=F(gf)\circ\eta_{g,X}$ for
    any $g\in G$ and $f\colon X\to Y$ a map in $\Cc$.
  \end{enumerate}
  Then, we will say that the family $\eta$ defines an action of $G$ on
  the functor $F$, or more succinctly, that $G$ acts on the functor
  $F$, or, following~\cite{MR2002h:55017}, that $F$ is a right
  $G$-functor. 
\end{definition}

\begin{definition}
  \label{definition:2}
  Let $F_{1},F_{2}\colon \Cc\to\Dd$ be two functors on which $G$ acts,
  by $\eta^{1},\eta^{2}$ respectively. A morphism of functors with
  $G$-action is a natural transformation $\epsilon\colon F_{1}\to
  F_{2}$ such that
  $\eta^{2}_{g,X}\circ\epsilon_{X}=\epsilon_{gX}\circ\eta^{1}_{g,X}$
  for all $g\in G$, $X\in\obj\Cc$. 
\end{definition}

As it is mentioned in \cite{MR2002h:55017} and \cite{myphdthesis}, the
usefulness of this concept lies on the fact that, when $G$ acts on
$F$, there is a natural action of $G$ on the simplicial set $\hocolim
F$, and in several other structures defined in terms of $F$. On the
other hand, it is often the case that we can derive a functor on which
$G$ acts by natural transformations from a $G$-object. We show some
examples. 

\begin{example}
  \label{example:1}
  Let $\Cc$ and $\Dd$ be $G$-categories, and $F\colon \Cc\to\Dd$ an
  equivariant functor. Then, for $D,D'\in\obj\Dd$, we have a category
  $\commac{D}{F}{D'}$ with objects
  $\obj(\commac{D}{F}{D'})=\setof{(u,C,v)}{C\in\obj\Cc,D\xrightarrow{u}
    FC\xrightarrow{v} D'}$, and a morphism $p\colon (u,C,v)\to
  (u',C',v')$ given by a $\Cc$-map $p\colon C\to C'$ such that
  $F(p)\circ u=u'$ and $v'\circ F(p)=v$.

  There is a functor $\Dd^{\op}\times \Dd\to\SCat$ defined on objects
  by $(D,D')\mapsto \commac{D}{F}{D'}$. If $(\phi,\psi)\colon
  (D,D')\to (E,E')$ is a morphism in $\Dd^{\op}\times \Dd$, the
  associated functor $\commac{D}{F}{D'}\to\commac{E}{F}{E'}$ sends
  $(u,C,v)$ to $(u\phi,C,\psi v)$.

  Then $\Dd^{\op}$ and $\Dd^{\op}\times \Dd$ have an obvious structure
  of $G$-categories, and there is an action of $G$ on the functor
  $\Dd^{\op}\times \Dd\to\SCat$ we just defined: for $g\in G$, set
  $\eta_{g,(D,D')}$ as the functor
  $\commac{D}{F}{D'}\to\commac{gD}{F}{gD'}$ given by $(u,C,v)\mapsto
  (gu,gC,gv)$.  Note, for example, that for $u\colon D\to F(C)$, we
  have that $gu\colon gD\to gF(C)=F(gC)$. 

  In this context, we can also define categories $\comma{D}$ and
  $\commaop{D}$ with the obvious objects and morphisms, and obtain
  functors $\Dd^{\op}\to\SCat$, $\Dd\to\SCat$ with a $G$-action. If
  $\nu\colon F_{1}\to F_{2}$ is an equivariant natural transformation
  (i.e.~a natural transformation such that $g\nu_{C}=\nu_{gC}$), then
  there is an induced morphism of right $G$-functors
  $\bar{\nu}\colon\comma[F_{1}]{-}\to\comma[F_{2}]{-}$, given by
  $\bar{\nu}_{D}\colon\comma[F_{1}]{D}\to\comma[F_{2}]{D}$,
  $(u,C)\mapsto (\nu_{C}u,C)$. 
\end{example}

\begin{example}
  \label{example:2}
  Again, let $\Cc$ and $\Dd$ be $G$-categories, and $F\colon
  \Cc\to\Dd$ an equivariant functor. There is a functor
  $\Cc^{\op}\times\Cc\to\Set$ defined on objects by
  $(X,Y)\mapsto\hom_{\Dd}(FX,FY)$ and on morphisms by
  $(\phi,\psi)\mapsto (f\mapsto F\phi\circ f\circ F\psi)$. It has a
  $G$-action defined by $\eta_{g,(X,Y)}\colon
  \hom_{\Dd}(FX,FY)\to\hom_{\Dd}(gFX,gFY)$, $f\mapsto gf$.

  Since any set $X$ can be considered as a simplicial set $Y$ such
  that $Y_{n}=X$ for all $n$ and all faces and degeneracies equal to
  the identity, we can as well consider the last functor as taking
  values in the category of simplicial sets. 
\end{example}

\begin{example}
  \label{example:5}
  Let $\Cc$ be a $G$-category and $F\colon \Cc\to\Dd$ a right
  $G$-functor with $G$-action given by $\eta$. Assume that $F$ has a
  colimit, that is, there is an object $\colim F$ in $\Dd$ and a
  collection of $\Dd$-maps $\{\rho_{X}\colon FX\to\colim
  F\}_{X\in\obj\Cc}$ that form a limiting cone from $F$ with base
  $\colim F$~(see for example \cite[p.~67]{categories}). Let $g\in G$.
  Then the natural transformation $F\to Fg$ induces a map $g\colon
  \colim F\to\colim Fg\cong\colim F$ such that
  $\rho_{gX}\circ\eta_{g,X}=g\circ\rho_{X}$ for all $X\in\obj\Cc$. It
  can be shown that the collection of maps $\{g\colon \colim
  F\to\colim F\}_{g\in G}$ give an structure of $G$-object on $\colim
  F$. Furthermore, if $Z$ is any $G$-object in $\Dd$ and there is a
  cone $\{\sigma_{X}\colon FX\to Z\}$ from~$F$ to $Z$ such that
  $\sigma_{gX}\circ\eta_{g,X}=g\circ\sigma_{X}$ for all $g\in G$ and
  all $X\in\obj\Cc$, then the map induced by the properties of the
  colimit $M\colon \colim F\to Z$ is in fact equivariant. 

  For example, if $\Cc$ is a discrete small $G$-category, then it can
  be identified with a $G$-set. A functor $F\colon \Cc\to\Dd$
  corresponds to a collection of $\Dd$-objects, indexed by the objects
  of $\Cc$. If $F$ is a right $G$-functor, then $\colim
  F=\coprod_{C\in\obj\Cc} F(C)$ is a $G$-object. 
  
  As a particular case, consider $H\le G$ a subgroup, and let $\Cc$ be
  the discrete small $G$-category with object set
  $G/\!/H=\{a_{1}H,a_{2}H,\ldots,a_{n}H\}$, that is, the set of left
  cosets of $H$ in $G$ with the usual action by left translation,
  where $a_{1}=1$. Let $Z$ be an $H$-simplicial set, and consider the
  constant functor $F\colon \Cc\to\SimpSets$ with value $Z$. We define
  a $G$-action $\eta$ on $F$ as follows: Let $\eta_{g,H}\colon F(H)\to
  F(gH)$ be defined as $z\mapsto hz$, where $g=a_{i}h$, $h\in H$; and
  then $\eta_{g,aH}(z)=\eta_{ga,H}(z)$. It is straightforward to check
  that this defines an action of $G$ on $F$, and so $\colim F$ is a
  $G$-simplicial set. This construction is usually known as the
  induced action from $H$ to $G$. We will denote $\colim F$ in this
  case as $\induc{Z}{H}{G}$.

  We also note that a morphism of right $G$-functors induces an
  equivariant map between the corresponding colimits of the functors. 
\end{example}

\begin{example}
  \label{example:6}
  In a similar way, if $Z\colon \Cc\times\Cc^{\op}\to \Dd$ is a right
  $G$-functor with a coend (see \cite[p.~226]{categories}) with
  defining maps $\alpha_{C}\colon Z(C,C)\to\coend Z$, then $\coend Z$
  becomes a $G$-object, with action satisfying
  $\alpha_{gC}\circ\eta_{g,(C,C)}=g\circ\alpha_{C}$. 

  For example, let $F\colon \Cc\to\SimpSets$, $T\colon
  \Cc^{\op}\to\SimpSets$ be functors, with actions of $G$ on both $F$
  and $T$, given by $\eta^{F},\eta^{T}$. Then $Z=F\times T$ is a right
  $G$-functor $\Cc\times\Cc^{\op}\to\SimpSets$. Its coend is a
  $G$-simplicial set denoted by $F\otimes_{\Cc} T$.

  As in the case of limits, a morphism of right $G$-functors induces
  an equivariant map between the corresponding coends. 
\end{example}

\begin{example}
  \label{example:4}
  Let $\Cc$ be a $G$-category and $Z\colon
  \Cc\times\Cc^{\op}\to\SimpSets$ a functor, with an action of $G$ on
  $Z$ given by $\eta$. We have a simplicial set $B(\Cc,Z)$, called the
  (simplicial) \emph{bar construction} (see \cite{MR86g:18010}), such
  that
  \begin{align}
    \label{eq:3}
    B(\Cc,Z)_{n} &=\coprod_{X_0\xrightarrow{\phi_1}
      X_1\xrightarrow{\phi_2}\cdots\xrightarrow{\phi_n} X_n\in N(\Cc)_{n}}
    Z(X_{0},X_{n})_{n}\\
    &=\setof{(\phi_{1},\ldots,\phi_{n};z)}{z\in Z(X_{0},X_{n})_{n}}
  \end{align}
  with boundaries and degeneracies given by:
  \begin{align}
    \label{eq:4}
    d^{i}(\phi_{1},\ldots,\phi_{n};z) &=
    \begin{cases}
      \(\phi_{2},\ldots,\phi_{n};d^{0}(Z(\phi_{1},1_{X_n})(z))\) &
      i=0,\\
      (\phi_{1},\ldots,\phi_{i+1}\phi_{i},\ldots,\phi_{n};d^{i}z)
      & 1\le i\le n-1,\\
      \(\phi_{1},\ldots,\phi_{n-1};d^{n}(Z(1_{X_0},\phi_{n})(z))\) & i=n
    \end{cases}\\
    s^{i}(\phi_{1},\ldots,\phi_{n};z) &=
    (\phi_{1},\ldots,\phi_{i},1_{X_i},\phi_{i+1},\ldots,\phi_{n};s^{i}z),\qquad
    0\le i\le n.
  \end{align}

  The action of $G$ on $Z$ gives a structure of $G$-simplicial set on
  $B(\Cc,Z)$, with action of $g\in G$ defined as:
  \begin{equation}
    \label{eq:5}
    g(\phi_{1},\ldots,\phi_{n};z)=
    \(g\phi_{1},\ldots,g\phi_{n};\eta_{g,(X_{0},X_n)}(z)\) 
  \end{equation}

  If the functor $Z$ is of the form $F\times T$ as in the previous
  example, then we denote $B(\Cc,Z)$ as $B(F,\Cc,T)$.
\end{example}

\begin{example}
  \label{sec:definition-examples-1}
  Let $F\colon\Cc\to\SCat$ be a functor. We define a category
  $\Cc\smallint F$ with objects the pairs $(X,a)$ with $X\in\obj\Cc$,
  $a\in\obj F(X)$. A map $(X,a)\to (Y,b)$ is given by a pair $(f,u)$
  such that $f\colon X\to Y$ is a map in $\Cc$ and $u\colon F(f)(a)\to
  b$ is a map in the category $F(Y)$.  The category $\Cc\smallint F$
  is called the \emph{Grothendieck Construction} on $F$
  (see~\cite{thomason}). 

  If $F\colon\Cc\to\SCat$ is a right $G$-functor, then $\Cc\smallint
  F$ is a small $G$-category with action on objects given by
  \begin{equation}
    \label{eq:6}
    g(X,a)=\(gX,\eta_{g,X}(a)\)
  \end{equation}
  and on maps by
  \begin{equation}
    \label{eq:8}
    g\((X,a)\xrightarrow{(f,u)}(Y,b)\)
    =\(gf,\eta_{g,Y}(u)\)
  \end{equation}  
\end{example}

% Our goal in the next section will be to generalize some of the
% theorems in~\cite{MR93e:55015} to this equivariant setting. 

We end this section by stating some basic and easily provable properties
of right $G$-functors.

\begin{proposition}
  \label{proposition:1}
  If $F\colon\Cc\to\Dd$ is functor with a $G$-action given by $\eta$
  and $X$ is an object in $\Cc$, then $FX$ is a $G_{X}$-object, where
  $G_{X}$ is the stabilizer of $X$ under the action of $G$ on
  $\obj\Cc$. The action is defined by the maps $\eta_{g,X}\colon FX\to
  FX$.
\end{proposition}

\begin{proposition}
  \label{proposition:2}
  {\normalfont ((2.3) from~\cite{MR2002h:55017})} Let
  $F\colon\Cc\to\Dd$ be a right $G$-functor, $S\colon \Cc'\to\Cc$ an
  equivariant functor, and $T\colon \Dd\to\Ee$ any functor. Then both
  $F\circ S$ and $T\circ F$ have induced structures of right
  $G$-functors. 
\end{proposition}

For example, for any $G$-category $\Cc$, we have a right $G$-functor
$N(\comma[\Cc]{-})\colon \Cc\to\SimpSets$. 

\section{The Homotopy Colimit}
\label{sec:homotopy-colimit}

% In the situation where $F\colon \Cc\times\Cc^{\op}\to\SimpSets$ is a
% functor, then for each $C'\in\obj\Cc$ one has a fuctor
% $F(\cdot,C')\colon \Cc\to\SimpSets$. As in~\cite[V.~3]{categories}, we
% obtain that both $B(T,\Cc,F)$ and $T\otimes_{\Cc} F$ are functors
% $\Cc^{\op}\to\SimpSets$.

Let $\Cc$ be a $G$-category and $Z\colon
\Cc\times\Cc^{\op}\to\SimpSets$ a right $G$-functor. We start by
noting the equivariant isomorphism:

\begin{equation}
  \label{eq:1}
  Z\otimes_{\Cc\times\Cc^{\op}}N(\commac{-}{\Cc}{-})\cong_{G} B(\Cc,Z),
\end{equation}
which can be proven by showing that $B(\Cc,Z)$ satisfies the
definition of coend of the functor $Z\times
N(\commac{-}{\Cc}{-})\colon (\Cc\times\Cc^{\op})\times
(\Cc\times\Cc^{\op})^{\op}\to\SimpSets$. In the case that $Z=F\times
T$ with $F\colon \Cc\to\SimpSets$, $T\colon \Cc^{\op}\to\SimpSets$ are
right $G$-functors, and using Fubini's theorem for
coends~\cite[p.~230]{categories}, this leads to 
\begin{equation}
  \label{eq:2}
  F\otimes_{\Cc} N(\commac{-}{\Cc}{-})\otimes_{\Cc} T\cong_{G} B(F,\Cc,T).
\end{equation}

Using that, we can prove that for right $G$-functors $F\colon
\Dd\to\SimpSets$, $T\colon \Cc^{\op}\times\Dd\to\SimpSets$ and
$U\colon \Dd^{\op}\to\SimpSets$, we have
\begin{equation}
  \label{eq:22}
  B(B(F,\Cc,T),\Dd,U)\cong_{G} B(F,\Cc,B(T,\Dd,U)),  
\end{equation}
whose non-equivariant version is~3.1.3 from~\cite{MR93e:55015}.

If $\Cc$ is any $G$-category, we will denote by $*$ the functor
$\Cc\to\SimpSets$ that is constant with value the simplicial set with
exactly one simplex in each dimension.  It is clearly has a structure
of right $G$-functor. 

\begin{definition}
  \label{sec:definition-examples-2}
  Let $F\colon \Cc\to\SimpSets$ a functor. Its \emph{homotopy colimit}
  $\hocolim_{\Cc} F$ is defined as $F \otimes_{\Cc} N(\comma[\Cc]{-})$.
\end{definition}

If $F$ is a right $G$-functor, then $Z=F\times N(\comma[\Cc]{-})$ has
a natural structure of right $G$-functor, so in this case
$\hocolim_{\Cc} F=\coend Z$ is a $G$-simplicial set. 

Note that the map of right $G$-functors $N(\comma[\Cc]{-})\to *$
induces an equivariant map
\begin{equation}
  \label{eq:9}
  \hocolim_{\Cc} F =F \otimes_{\Cc} N(\comma[\Cc]{-})\to F\otimes_{\Cc}
  *=\colim F,
\end{equation}
and the map of right $G$-functors $F\to *$ induces an equivariant map
\begin{equation}
  \label{eq:10}
  \hocolim_{\Cc} F =F \otimes_{\Cc} N(\comma[\Cc]{-})\to
  N(\comma[\Cc]{-})\otimes_{\Cc}*
  =N(\Cc)
\end{equation}

One also can prove the isomorphism of right $G$-functors:
\begin{equation}
  \label{eq:11}
  N(\commac{-}{\Cc}{-})\otimes_{\Cc^{\op}} *\cong N(\comma[\Cc]{-})
\end{equation}
which leads to the equivariant isomorphism:
\begin{equation}
  \label{eq:12}
  B(F,\Cc,*)\cong_{G}\hocolim_{\Cc} F.
\end{equation}

Finally, we note that just by categorical arguments, one obtains:

\begin{proposition}
  \label{proposition:3}
  Let $S\colon \Dd\to\Cc$ be an equivariant functor between
  $G$-categories, and let $F\colon \Cc\to\SimpSets$ a right
  $G$-functor. Then, with the induced right $G$-functor structure on
  $F\circ S$, we have:
  \begin{enumerate}
  \item \label{item:1} $\hom_{\Cc}(C,S-)\otimes_{\Dd}
    N(\comma[\Dd]{-})\cong B(\hom_{\Cc}(C,S-),\Dd,*)\cong
    N(\comma[S]{C})$ as right $G$-functors on the argument $C$. 
  \item \label{item:2} $(F\circ S)(D)\cong
    F\otimes_{\Cc}\hom_{\Cc}(-,SD)\cong B(F,\Cc,\hom_{\Cc}(-,SD))$, as
    right $G$-functors on $D$. 
  \end{enumerate}
\end{proposition}

As a consequence of this proposition, if we take $S=1_{\Cc}$ to be the
identity functor, we obtain that 
\begin{equation}
  \label{eq:28}
  B(\hom_{\Cc}(C,-),\Cc,*)\cong_{G_C} N(\comma[\Cc]{C})\simeq_{G_C} *,
\end{equation}
for all $C\in\obj\Cc$, since $\comma[\Cc]{C}$ has an initial object
$1_{C}\colon C\to C$ fixed by $G_{C}$
(\cite[(4.3)]{equivalence:simplicial}). 

\section{The Homotopy Invariance Theorem}
\label{sec:2}

The proofs of the theorems of the next section are based on this
important theorem. The reader may refer to~\cite{tomdieck} for the
properties of induced topological spaces. 

\begin{theorem}
  \label{theorem:1}
  Let $Z,Z'\colon \Cc\times\Cc^{\op}\to\SimpSets$ two right
  $G$-functors. Let $\epsilon\colon Z\to Z'$ be a map of right
  $G$-functors such that $\epsilon_{(X,Y)}\colon Z(X,Y)\to Z'(X,Y)$ is
  a $G_{(X,Y)}$-homotopy equivalence for all
  $X\in\obj\Cc,Y\in\obj\Cc^{\op}$. Then the map $\bar{\epsilon}$
  induced by $\epsilon$:
  \begin{equation}
    \label{eq:13}
    \bar{\epsilon}\colon B(\Cc,Z)\to B(\Cc,Z')
  \end{equation}
  is a $G$-homotopy equivalence.
\end{theorem}

\begin{proof}
  From~\cite{MR86g:18010}, we know that $B(\Cc,Z)$ is the diagonal of
  a bisimplicial set ${\tilde B}(\Cc,Z)$ with $(m,n)$-simplices the
  set
  \begin{equation}
    \label{eq:14}
    \coprod_{X_0\xrightarrow{\phi_1}
      X_1\xrightarrow{\phi_2}\cdots\xrightarrow{\phi_m} X_m}
    Z(X_{0},X_{m})_{n}.
  \end{equation}
  From the examples, we know that this coproduct has an action of $G$
  given by:
  \begin{equation}
    \label{eq:16}
    g(\phi_{1},\ldots,\phi_{m};z)=
    \(g\phi_{1},\ldots,g\phi_{m};\eta_{g,(X_{0},X_m)}(z)\),
  \end{equation}
  and this makes ${\tilde B}(\Cc,Z)$ a bisimplicial $G$-set.
  We have that $\epsilon$ induces a map ${\tilde \epsilon}\colon {\tilde
    B}(\Cc,Z)\to{\tilde B}(\Cc,Z')$, sending
  \begin{equation}
    \label{eq:15}
    (\phi_{1},\ldots,\phi_{m};z)\mapsto
    \(\phi_{1},\ldots,\phi_{m};\epsilon_{(X_0,X_m)}(z)\),
  \end{equation}
  The map ${\tilde \epsilon}$ is equivariant, and so if we define
  $\bar{\epsilon}$ as $\diag{\tilde \epsilon}$, then $\bar{\epsilon}$
  is equivariant as well.

  Let us denote $X_0\xrightarrow{\phi_1}
  X_1\xrightarrow{\phi_2}\cdots\xrightarrow{\phi_m} X_m\in N(\Cc)_{m}$
  by $\bar{X}$. According to
  Theorem~(3.8) from~\cite{equivalence:simplicial}, in order to prove
  that $\bar{\epsilon}$ is a $G$-homotopy equivalence, it is
  sufficient to prove that
  \begin{equation}
    \label{eq:17}
    {\tilde \epsilon}_{m,-}\colon
    \coprod_{\bar{X}\in N(\Cc)_{m}} Z(X_{0},X_{m})
    \to
    \coprod_{\bar{X}\in N(\Cc)_{m}} Z'(X_{0},X_{m})
  \end{equation}
  is a $G$-homotopy equivalence for all $m$. Taking geometric
  realization on both sides of~\eqref{eq:17}, since geometric
  realization commutes with coproducts, we obtain:
  \begin{equation}
    \label{eq:18}
    |{\tilde \epsilon}_{m,-}|\colon
    \coprod_{\bar{X}\in N(\Cc)_{m}} |Z(X_{0},X_{m})|
    \to
    \coprod_{\bar{X}\in N(\Cc)_{m}} |Z'(X_{0},X_{m})|
  \end{equation}
  Let $E_m$ be a set of representatives for the orbits of the action
  of $G$ on $N(\Cc)_m$. Then the map in~\eqref{eq:18} can be written
  as:
  \begin{equation}
    \label{eq:19}
    \induc{|{\tilde \epsilon}_{m,-}|}{G_{\bar{Y}}}{G}\colon 
    \coprod_{\bar{Y}\in E_m} \induc{|Z(X_{0},X_{m})|}{G_{\bar{Y}}}{G}\to
    \coprod_{\bar{Y}\in E_m} \induc{|Z'(X_{0},X_{m})|}{G_{\bar{Y}}}{G}
  \end{equation}
  Since by hypothesis, each $\epsilon_{(X_{0},X_{m})}$ is a
  $G_{(X_{0},X_{m})}$-homotopy equivalence, given that $G_{\bar{Y}}\le
  G_{(X_{0},X_{m})}$, they are also $G_{\bar{Y}}$-homotopy
  equivalences, and so each map $
  \induc{|Z(X_{0},X_{m})|}{G_{\bar{Y}}}{G}\to
  \induc{|Z'(X_{0},X_{m})|}{G_{\bar{Y}}}{G}$ is a $G$-homotopy
  equivalence. Therefore the map in~\eqref{eq:19} is a coproduct of
  $G$-homotopy equivalences, hence a $G$-homotopy equivalence, as we
  wanted to prove. 
\end{proof}

\section{Further Theorems}
\label{sec:3}

\begin{theorem}
  \label{theorem:2}
  {\normalfont (Equivariant Homotopy Invariance Of The Homotopy
    Colimit)}. Let $F,F'\colon \Cc\to\SimpSets$ right $G$-functors,
  and $\epsilon\colon F\to F'$ a map of right $G$-functors such that
  each $\epsilon_{X}\colon FX\to F'X$ is a $G_{X}$-homotopy
  equivalence. Then the induced map $\bar{\epsilon}\colon
  \hocolim_{\Cc} F\to\hocolim_{\Cc} F'$ is a $G$-homotopy equivalence. 
\end{theorem}

\begin{proof}
  Straightforward from Theorem~\ref{theorem:1}, since the homotopy
  colimit is a special case of a bar construction.
\end{proof}

\begin{theorem}
  \label{theorem:3}
  {\normalfont (Reduction Theorem)} Let $S\colon \Dd\to\Cc$ be an
  equivariant functor between $G$-categories, and let $F\colon
  \Cc\to\SimpSets$ a right $G$-functor. Then we have the equivariant
  isomorphism.
  \begin{equation}
    \label{eq:7}
    \hocolim_{\Dd} F\circ S\cong_{G} F\otimes_{\Cc}N(\comma[S]{-})
  \end{equation}
\end{theorem}

\begin{proof}
  \begin{align*}
    \hocolim_{\Dd} F\circ S &= (F\circ S)\otimes_{\Dd}
    N(\comma[\Dd]{-})
    &&\text{Definition of $\hocolim$}\\
    &\cong_{G} (F\otimes_{\Cc}\hom_{\Cc}(-,SD))\otimes_{\Dd}
    N(\comma[\Dd]{-})
    &&\text{Proposition~\ref{proposition:3}.\ref{item:2}}\\
    &\cong_{G} F\otimes_{\Cc}(\hom_{\Cc}(C,S-)\otimes_{\Dd}
    N(\comma[\Dd]{-}))
    &&\text{Fubini's theorem}\\
    &\cong_{G} F\otimes_{\Cc}N(\comma[S]{-})
    &&\text{Proposition~\ref{proposition:3}.\ref{item:1}}\qedhere
  \end{align*}
\end{proof}

In~\cite[(2.6)]{MR2002h:55017}, this result is given as a homotopy
equivalence.  However, as noted in~\cite[4.4]{MR93e:55015}, this is
even an isomorphism, which in this case is equivariant.

\begin{theorem}
  \label{theorem:5}
  {\normalfont (Cofinality Theorem)} Let $S\colon \Dd\to\Cc$ be an
  equivariant functor between $G$-categories, and let $F\colon
  \Cc\to\SimpSets$ a right $G$-functor. Consider the induced right
  $G$-functor structure on $F\circ S$. If $N(\comma[S]{C})$ is
  $G_{C}$-contractible for all objects $C$ in $\Cc$, then
  $\hocolim_{\Dd}F\circ S\simeq_{G}\hocolim_{\Cc} F$. 
\end{theorem}

\begin{proof}
  \begin{align*}
    \hocolim_{\Dd}F\circ S &= B(F\circ S,\Dd,*)
    &&\text{Equation~\ref{eq:12}}\\
    &=B(B(F,\Cc,\hom_{\Cc}(-,SD)),\Dd,*)
    &&\text{Proposition~\ref{proposition:3}.\ref{item:2}}\\
    &\cong_{G} B(F,\Cc,B(\hom_{\Cc}(C,S-),\Dd,*))
    &&\text{Equation~\ref{eq:22}}\\
    &\cong_{G} B(F,\Cc,N(\comma[S]{-}))
    &&\text{Proposition~\ref{proposition:3}.\ref{item:1}}\\
    &\simeq_{G} B(F,\Cc,*)=\hocolim_{\Cc}F
    &&\text{Hypothesis}\qedhere
  \end{align*}
\end{proof}

We note that the hypothesis about $G_{C}$-contractibility of the fiber
$N(\comma[S]{C})$ allows us to conclude the $G$-homotopy. Compare
with~\cite[(2.7)]{MR2002h:55017}, where this result is given as a
homotopy equivalence not necessarily equivariant. 

\begin{theorem}
  \label{theorem:4}
  {\normalfont (Homotopy Pushdown Theorem)} Let $S\colon \Dd\to\Cc$ be
  an equivariant functor and $F\colon \Dd\to\SimpSets$ a right
  $G$-functor. Let $S_{h_*}(F)\colon \Cc\to\SimpSets$ the functor
  given by $C\mapsto B(F,\Dd,\hom_{\Cc}(S-,C))$. Then $S_{h_*}(F)$ is
  a right $G$-functor and $\hocolim_{\Cc}
  S_{h_*}(F)\simeq_{G}\hocolim_{\Dd}F$. 
\end{theorem}

\begin{proof}
  \begin{align*}
    \hocolim_{\Cc}S_{h_*}(F) &=B(B(F,\Dd,\hom_{\Cc}(S-,C)),\Cc,*)
    &&\text{Definition}\\
    &\cong_{G} B(F,\Dd,B(\hom_{\Cc}(SD,-),\Cc,*))
    &&\text{Equation~\ref{eq:22}}\\
    &\simeq_{G} B(F,\Dd,*)=\hocolim_{\Dd} F
    &&\text{Equation~\ref{eq:28}}\qedhere
  \end{align*}
\end{proof}

Note that we also used the equivariant homotopy invariance
(Theorem~\ref{theorem:1}) of the bar construction in the last step.
Hence in~\cite[(2.5)]{MR2002h:55017} we do have a $G$-homotopy
equivalence. 

\bibliography{mybiblio}
\bibliographystyle{amsplain}

\end{document}